\documentclass[12pt]{article}
\usepackage{amsmath,amsthm}
\usepackage{amssymb,latexsym}
\usepackage{enumerate}
\usepackage[english]{babel}
\usepackage{graphicx}

\def\Q{{\mathbb Q}}
\def\Z{{\mathbb Z}}

\newtheorem{lemma}{Lemma}
\newtheorem{theorem}[lemma]{Theorem}

\newtheorem{remark}[lemma]{Remark}

\title{
Integral bases and monogenity of pure fields\\
}
\author{
Istv\'{a}n Ga\'{a}l\thanks{
        Research supported in part by K115479 from the
        Hungarian National Foundation for Scientific Research
                         },\; 
and L\'aszl\'o Remete
\\ \\
University of Debrecen, Mathematical Institute \\
H--4002 Debrecen Pf.400., Hungary \\
e--mail: gaal.istvan@unideb.hu, remete.laszlo@science.unideb.hu
}

\begin{document}

\maketitle
\thispagestyle{empty}

\renewcommand{\thefootnote}{}



\renewcommand{\thefootnote}{\arabic{footnote}}
\setcounter{footnote}{0}

\begin{abstract}
Let $m$ be a square-free integer ($m\neq 0,\pm 1$).
We show that the structure of the integral bases of the fields
$K=\Q(\sqrt[n]{m})$ are periodic in $m$.
For $3\leq n\leq 9$ we show that the period length is $n^2$.
We explicitly describe the integral bases,
and for $n=3,4,5,6,8$ we explicitly calculate the index forms of $K$.
This enables us in many cases to characterize
the monogenity of these fields. 
Using the explicit form of the index forms yields a new
technic that enables us to derive new results on
monogenity and to get several former results as easy consequences.
For $n=4,6,8$ we give an almost complete characterization 
of the monogenity of pure fields.

\end{abstract}

\section{Introduction}

Let $m$ be a square-free integer ($m\neq 0,\pm 1$) 
and $n\geq 2$ a positive integer. 
There is an extensive literature of pure fields of type $K=\Q(\sqrt[n]{m})$. 
(Describing the following results on pure fields we use some basic 
concepts on monogenity and power integral bases
that are detailed in Section \ref{iiii}.)

B.K.Spearman and K.S.Williams \cite{sp} gave an explicit formula 
for the integral basis of pure cubic fields. 
B.K.Spearman, Y.Qiduan and J.Yoo \cite{sqy} showed that
if $i$ is a cubefree positive integer then there exist infinitely
many pure cubic fields with minimal index equal to $i$.
I.Ga\'al and T.Szab\'o \cite{gsz} studied the behaviour of the minimal
indices of pure cubic fields in terms of the discriminant.
L. El Fadil \cite{fad} gave conditions for the existence of 
power integral bases of pure cubic fields
in terms of the index form equation.

T.Funakura \cite{fun} studied the integral basis in pure quartic fields.
I.Ga\'al and L.Remete \cite{gr} calculated elements of 
index 1 (with coefficients $<10^{1000}$) in pure quartic fields 
$K=\Q(\sqrt[4]{m})$ for $1<m<10^7$, $m\equiv 2,3 \; (\bmod \ 4)$.

S.Ahmad, T.Nakahara and S.M.Husnine \cite{n2} showed that if
$m\equiv 1 \; (\bmod \ 4)$,  $m\not\equiv \pm 1 \; (\bmod \ 9)$ then 
$\Q(\sqrt[6]{m})$ is not monogenic. On the other hand  \cite{n1}, if
$m\equiv 2,3 \; (\bmod \ 4)$,  $m\not\equiv \pm 1 \; (\bmod \ 9)$ then 
$\Q(\sqrt[6]{m})$ is monogenic.

A.Hameed and T.Nakahara \cite{n3} constructed integral bases of pure octic
fields $\Q(\sqrt[8]{m})$. They proved \cite{n8} that
if $m\equiv 1\; (\bmod \ 4)$ then $\Q(\sqrt[8]{m})$
is not monogenic. On the other hand 
A.Hameed, T.Nakahara, S.M.Husnine and S.Ahmad \cite{nnn} proved that 
if $m\equiv 2,3\; (\bmod \ 4)$  then $\Q(\sqrt[8]{m})$
is monogenic.

A.Hameed, T.Nakahara, S.M.Husnine and S.Ahmad \cite{nnn} showed that
if $m\equiv 2,3\; (\bmod \ 4)$  then $\Q(\sqrt[2^n]{m})$ is monogenic,
this involves the pure quartic and pure octic fields, as well.
Moreover, they showed \cite{nnn} that
if all the prime factors of $n$ divide $m$ 
then $\Q(\sqrt[n]{m})$ is monogenic.

\vspace{1cm}

Our purpose is for $3\leq n\leq 9$ to 
give a general characterization of the integral basis 
of $K=\Q(\sqrt[n]{m})$. 
We prove that the integral bases 
of $K=\Q(\sqrt[n]{m})$ are periodic in $m$.
For $3\leq n\leq 9$ the period length is $n^2$.

The knowledge of the integral bases makes 
possible also to compute the
sporadic results on the monogenity of these fields.
Our method applying the explicit form of the index forms yields
a new technics that enables us to obtain new results on the
monogenity of these fields and to obtain several former results
as easy consequences.

In our Theorems \ref{tquartic}, \ref{tsextic}, \ref{toctic}
we give an almost complete characterization of the
monogenity of pure quartic, sextic and octic fields, respectively.
The cubic case is well-known and easy, much less is known
about the quintic, septic and nonic cases.

\section{Basic concepts about the monogenity of number fields}
\label{iiii}

We recall those concepts \cite{gaal} that we use throughout.
Let $\alpha$ be a primitive integral element of the number field $K$
(that is $K=\Q(\alpha)$) of degree $n$ with ring of integers $\Z_K$. 
The {\bf index} of $\alpha$ is
\[
I(\alpha)=(\Z_K^{+}:\Z[\alpha]^{+})=\sqrt{\left|\frac{D(\alpha)}{D_K}\right|}
=\frac{1}{\sqrt{|D_K|}}\prod_{1\leq i<j\leq n}\left|\alpha^{(i)}-\alpha^{(j)}\right| \;\; ,
\]
where $D_K$ is the discriminant of $K$ and $\alpha^{(i)}$ denote
the conjugates of $\alpha$. The minimal index of $K$ is 
\[
i_K=\min I(\alpha)
\]
where $\alpha$ runs through the primitive integral elements of $K$.

If $B=(b_1=1,b_2,\ldots,b_n)$ is an
integral basis of $K$, then the {\bf index form}
corresponding to this integral basis is
\[
I(x_2,\ldots,x_n)=\frac{1}{\sqrt{|D_K|}}
\prod_{1\leq i<j\leq n}
\left(
(b_2^{(i)}-b_2^{(j)})x_2+\ldots+(b_n^{(i)}-b_n^{(j)})x_n
\right)
\]
(where $b_j^{(i)}$ denote the conjugates of $b_j$)
which is a homogeneous polynomial with integral coefficients.
For the integral element 
\[
\alpha=x_1+b_2x_2+\ldots+b_nx_n
\]
we have
\[
I(\alpha)=|I(x_2,\ldots,x_n)|
\]
independently of $x_1$.
$\alpha$ generates a {\bf power integral basis} $(1,\alpha,\ldots,\alpha^{n-1}$)
if and only if $I(\alpha)=1$ that is $(x_2,\ldots,x_n)\in \Z^{n-1}$ is 
a solution of the
{\bf index form equation}
\begin{equation}
I(x_2,\ldots,x_n)=\pm 1 \;\;\; {\rm in} \;\;\; (x_2,\ldots,x_n)\in \Z^{n-1}.
\label{iiixxx}
\end{equation}
In this case
\[
\Z_K=\Z[\alpha]
\] 
and $K$ is called {\bf monogenic}.

\section{Basic results}

Throughout we assume that $m$ is a square-free integer
with $m\neq 0,\pm 1$ and $n>2$ an integer.
Let $K=\Q(\sqrt[n]{m})$ and $\vartheta=\sqrt[n]{m}$.

\vspace{1cm}

Our first theorem is on the prime divisors of
the denominators of the integral basis elements:

\begin{theorem}
\label{th1}
If $(1,\vartheta,\ldots,\vartheta^{n-1})$ is not an integral basis
in $K$, then for any element 
\begin{equation}
\alpha=\frac{a_0+a_1\vartheta+\ldots+a_{n-1}\vartheta^{n-1}}{q}
\label{pp}
\end{equation}
of the integral basis 
(with $a_0,\ldots,a_{n-1},q\in\Z$, $q> 0$)
the denominator $q$ can only be divisible
by primes dividing $n$, the prime factors of $q$
do not divide $m$.
\end{theorem}

\noindent
{\bf Proof} \\
\noindent
The discriminant of $\vartheta=\sqrt[n]{m}$ is $\pm n^nm^{n-1}$.
If $(1,\vartheta,\ldots,\vartheta^{n-1})$ is not an integral basis
in $\Q(\vartheta)$, then
there must be a positive integer $q$ dividing $n^nm^{n-1}$ and an element 
$\alpha$ of type (\ref{pp})
such that $\alpha$ is an algebraic integer and an element of 
$(1,\vartheta,\ldots,\vartheta^{n-1})$ can be replaced by $\alpha$
to get a basis with smaller discriminant.

Let $p$ be a prime divisor of $q$. Then obviously 
\begin{equation}
\alpha'=\frac{q}{p}\alpha=\frac{e_0+e_1\vartheta+\ldots+e_{n-1}\vartheta^{n-1}}{p}
\label{ppv}
\end{equation}
is also an algebraic integer. 
We can also assume that $0\leq e_i<p \; (0\leq i\leq n-1)$ 
by taking each $e_i$ modulo $p$.

We show that $p$ is a divisor of $n$. 

Assume on the contrary that $p|m$.  The element
\[
\alpha'\vartheta=
\frac{e_0\vartheta+e_1\vartheta^2+\ldots+e_{n-2}\vartheta^{n-1}+e_{n-1}m}{p}
\]
is obviously an algebraic integer. By $p|m$ the element
\[
\frac{e_0\vartheta+e_1\vartheta^2+\ldots+e_{n-2}\vartheta^{n-1}}{p}
\]
is also an algebraic integer. We proceed by multiplying this element
by $\vartheta$ and omitting the analogous integral part. Finally we
obtain that 
\[
\varrho=\frac{e_0\vartheta^{n-1}}{p}
\]
is an algebraic integer.
The element $\varrho$ is the root of the polynomial
\[
f_{\varrho}(x)=px^n-e_0^n\left(\frac{m}{p}\right)^{n-1}.
\]
This polynomial is irreducible over $\Q$
if and only if its reciprocal polynomial
\[
f_{1/\varrho}(x)=e_0^n\left(\frac{m}{p}\right)^{n-1}x^n-p.
\]
is irreducible. 
Here $m/p$ is an integer, not divisible by $p$ because $m$
is square-free. $e_0$ is also not divisible by $p$ (otherwise
we did not have $e_0$ in (\ref{ppv}) and we had the same result 
with the first non-zero $e_i$). Hence $f_{1/\varrho}(x)$ is 
an Eisentein polynomial, therefore $f_{\varrho}(x)$ is irreducible.
Then $f_{\varrho}(x)$ is the defining polynomial of $\varrho$. 
This contradicts to $\varrho$ being an algebraic integer. $\Box$

\vspace{1cm}

{\bf Remark.} Theorem \ref{th1} implies Theorem 3.1. of 
A.Hameed, T.Nakahara, S.M.Husnine and S.Ahmad \cite{nnn}:
if all the prime factors of $n$ divide $m$ then $\Q(\sqrt[n]{m})$ is monogenic.

\vspace{1cm}

Next we show that the integral bases of $K=\Q(\sqrt[n]{m})$ are periodic.
First we prove this statement with a period 
length much larger than $n^2$ but
this result is valid for any $n$.

\begin{theorem}
\label{th2}
Let $n=p_1^{h_1}\ldots p_k^{h_k}$ and $n_0=p_1^{[nh_1/2]}\ldots p_k^{[nh_k/2]}$
where $[x]$ denotes the lower integer part of $x$.
Let $\vartheta=\sqrt[n]{m}$ and $\gamma=\sqrt[n]{m+n_0^n}$. Then
the structure of the integral bases of
the fields $\Q(\vartheta)$ and $\Q(\gamma)$ is the same
in terms of $\vartheta$ and $\gamma$, respectively.
\end{theorem}

\vspace{0.5 cm}

\noindent
{\em Remark.}
Under the "same structure" we mean that if the integral basis of
$\Q(\vartheta)$ has an element 
\[
\frac{a_0+a_1\vartheta+\ldots+a_{n-1}\vartheta^{n-1}}{q},
\]
then the integral basis of $\Q(\gamma)$ has an element 
\[
\frac{a_0+a_1\gamma+\ldots+a_{n-1}\gamma^{n-1}}{q}
\]
and vice versa.

\vspace{0.5 cm}

\noindent
{\bf Proof} \\
\noindent
Assume that $(1,\vartheta,\ldots,\vartheta^{n-1})$ is not an integral basis
in $\Q(\vartheta)$. 
Then there must be integer elements of type
\begin{equation}
\alpha=\frac{a_0+a_1\vartheta+\ldots+a_{n-1}\vartheta^{n-1}}{q}
\label{mm}
\end{equation}
which can replace elements of $(1,\vartheta,\ldots,\vartheta^{n-1})$
to obtain an integral basis. 
We show that the existence of analogous algebraic integers of type (\ref{mm}) is
equivalent in the fields generated by 
$\vartheta=\sqrt[n]{m}$ and by $\gamma=\sqrt[n]{m+n_0^n}$.

If we replace an element of the basis $(1,\vartheta,\ldots,\vartheta^{n-1})$
by $\alpha$ of (\ref{mm}) then the discriminant of the basis
decreases by a factor $q^2$. Hence $q^2$ divides $\pm n^nm^{n-1}$
(the discriminant of $\vartheta=\sqrt[n]{m}$ is $\pm n^nm^{n-1}$).
By Theorem \ref{th1} the prime divisors $p$ of $q$ do not divide $m$.
Hence $q^2|n^n$, which implies that $q$ divides $n_0$.

Denote the
conjugates of $\alpha$ by $\alpha^{(j)},j=1,\ldots,n$. 
The defining polynomial of $\alpha$ is
\[
\prod_{j=1}^n (x-\alpha^{(j)})=
\frac{1}{q^n}\prod_{j=1}^n 
(qx-a_0-a_1\vartheta^{(j)}-\ldots-a_{n-1}(\vartheta^{(j)})^{n-1}).
\]
The product is a symmetrical polynomial of $\vartheta^{(1)},\ldots,\vartheta^{(n)}$,
hence its coefficients can be expressed as polynomials (with integer
coefficients) of the defining polynomial of $\vartheta$, that is $x^n-m$.
Hence there exist polynomials $P_0,\ldots,P_{n-1}\in\Z[x]$ such that
\[
\prod_{j=1}^n (x-\alpha^{(j)})=
\frac{1}{q^n} ((qx)^n+P_{n-1}(m)(qx)^{n-1}+\ldots +P_1(m)(qx)+P_0(m)).
\]
Therefore the element $\alpha$ is an algebraic integer if and only if
$q^n|q^j P_j(m)$ that is 
\begin{equation}
q^{n-j}|P_j(m) \;\;\; (j=0,1,\ldots ,n-1).
\label{divp1}
\end{equation}

Replace now $m$ by $m'=m+n_0^n$ and consider integral bases
in the field $\Q(\gamma)=\Q(\sqrt[n]{m+n_0^n})$. 
In this field an element of type (\ref{mm}),
that is 
\begin{equation}
\delta=\frac{a_0+a_1\gamma+\ldots+a_{n-1}\gamma^{n-1}}{q}
\label{nn}
\end{equation}
is an algebraic integer if and only if
\begin{equation}
q^{n-j}|P_j(m+n_0^n)\;\; (j=0,1\ldots,n-1)
\label{divp2}
\end{equation}
with the same polynomials $P_j$. By $q|n_0$ the conditions
(\ref{divp1}) are equivalent to (\ref{divp2}). Therefore 
$\Q(\vartheta)$ and $\Q(\gamma)$ contains the same type of integer
elements. Elements of that type are linearly independent
in the first case if and only if they are linearly independent in the
second case. Therefore  $\Q(\vartheta)$ and $\Q(\gamma)$
admits the same type of integral bases. $\Box$

\vspace{1cm}

\section{The structure of the integral bases is \\
periodic in $m$ modulo $n^2$ for $3\leq n\leq 9$}

Theorem \ref{th2} implies that the integral bases of $K=\Q(\sqrt[n]{m})$
are periodic modulo $n_0^n$.
This number is of magnitude $n^{n^2/2}$. 
For small values of $n$ we have a much sharper assertion.

\begin{theorem}
For $3\leq n\leq 9$ the integral bases of $\Q(\sqrt[n]{m})$
are periodic in $m$ modulo $n^2$.
\label{th3}
\end{theorem}

\noindent
{\bf Proof} \\
\noindent
For $n=3,4,5$ the $n_0^n$ is 27, 655536, 9765625, respectively.
Calculating the integral bases of $\Q(\sqrt[n]{m})$ 
for square-free $m$ up to $n_0^n$ 
it is easily seen that the structure of the integral bases of $\Q(\sqrt[n]{m})$
are periodic modulo $n^2$. One can easily detect a few types of integral bases 
that are repeated for square-free values of $m$, $m+n^2$, $m+2n^2$ etc.

\vspace{1cm}

Let now $n>5$. Then $n_0^n$ is far too large for the calculations
described above.
However for $n=6,7,8,9$ we managed to prove the same assertion.

Let $1<r<n^2$. 
If $r$ is square-free, then set $r'=r$.
If $r$ has a common square factor with $n$,
then none of $r+kn^2$ is square-free, we omit $r$. 
If $r$ has no common square factor with $n$ but contains another square factor,
then we set
$r'=r+n^2$ or $r'=r+2n^2$ etc. which is already square-free.

Let $\vartheta=\sqrt[n]{r'}$,
calculate the integral bases of $\Q(\sqrt[n]{r'})$ and
denote the basis elements by 
$(b_1=1,b_2,\ldots,b_n)$, where $b_j$ is of the form
\[
b_j=\frac{a_{j0}+a_{j1}\vartheta+\ldots+a_{j,n-1}\vartheta^{n-1}}{q}
\]
(with $a_{j0},a_{j1},\ldots, a_{j,n-1}\in\Z$ and 
with a non-zero denominator $q$ the prime factors of which divide $n$).

Let $m=r+kn^2$ be a square-free integer, $\gamma=\sqrt[n]{m}$ and 
\[
b'_j=\frac{a_{j0}+a_{j1}\gamma+\ldots+a_{j,n-1}\gamma^{n-1}}{q}.
\]

\noindent
We wonder \\
I. if the analogues of the elements $b_j$, that is the elements $b'_j$ 
remain algebraic integer for any square-free $m=r+kn^2$, further \\
II. if for some square-free 
$m=r+kn^2$ some of the basis elements $(b_1'=1,b_2',\ldots,b_n')$
can be replaced by an integral element of type
\begin{equation}
d=\frac{e_1 b_1'+\ldots+e_n b_n'}{p}
\label{di}
\end{equation}
(where $0\leq e_1,\ldots,e_n\leq p-1$ and $p$ is a prime divisor of $n$)
to obtain a basis with smaller discriminant.

\vspace{1cm}

\noindent
I. The defining polynomial of $e_1 b_1'+\ldots+e_n b_n'$
is 
\[
G(x)=\prod_{j=1}^n (x-e_1 b_1'^{(j)}-\ldots-e_n b_n'^{(j)}).
\]
This polynomial is symmetrical in the conjugates of 
$\gamma=\sqrt[n]{m}$, hence its coefficients will be
polynomials in $m$:
\[
G(x)=x^n+G_{n-1}(m)x^{n-1}+\ldots +G_1(m)x+G_0(m).
\]
Since there are denominators in the $b_i'$, the polynomials
$G_j$ (depending also on $e_1,\ldots,e_n$)
are not necessarily of integer coefficients.
Let us substitute $m=r+kn^2$. We obtain
\[
G(x)=x^n+H_{n-1}(k)x^{n-1}+\ldots +H_1(k)x+H_0(k).
\]
For all possible residues $r$ we have explicitly calculated 
these polynomials $H_j(k)$ which also depend on $e_1,\ldots,e_n$.
In all cases we found that these are polynomials in
$k,e_1,\ldots,e_n$ with integer coefficients.
Substituting $e_i=1$ and $e_j=0,j=1,\ldots,n,\; j\neq i$
this implies that $b_i'$ is integer for any 
square-free $m=r+kn^2$ ($1\leq i\leq n$).

\vspace{1cm}

\noindent
II.Consider now the defining polynomial of $d$ (\ref{di}), that is
\[
P(x)=\frac{1}{p^n}G(px)=\frac{1}{p^n}
\left(  (px)^n+H_{n-1}(k)(px)^{n-1}+\ldots +H_1(k)(px)+H_0(k) \right).
\]
This polynomial has integer coefficients if and only if
$p^n$ divides $p^jH_j(k)$ that is
\begin{equation}
p^{n-j}| H_j(k)=H_j(k,e_1,\ldots,e_n) \;\; (0\leq j\leq n-1).
\label{test}
\end{equation}
Let now $(e_1,\ldots,e_n)\in \Z^n$ be an arbitrary given fixed tuple. 
Obviously by (\ref{test}) the validity of the statement if $d$ is integral or not,
depends only on the behaviour of $k$ modulo $p^n$ and not
on the value of $k$. This allows us to test the fields  
$\Q(\sqrt[n]{m})$ for square-free $m=r+kn^2$ where 
$k$ runs though all
residue classes modulo $p^n$. These fields were tested directly,
calculating their integral bases. We found that in all cases
the fields had the same structure of integral basis in terms 
of $\gamma=\sqrt[n]{m}$ like the field $\Q(\sqrt[n]{r'})$
in terms of $\vartheta=\sqrt[n]{r'}$.
This proves our assertion. $\Box$

\vspace{1cm}

\noindent
{\bf Remark} The test described at the end of the above proof
required to calculate the integral bases of \\
$24(2^6+3^6)=18239$ sextic fields, \\
$48\cdot 7^7=39530064$ septic fields,\\
$48\cdot 2^8=12288$ octic fields and\\
$72\cdot 3^9=1417176$ nonic fields.

\vspace{1cm}

\section{Integral bases and monogenity of pure fields}
\label{ppp}

In this section we give a list of our results on 
the fields $K=\Q(\sqrt[n]{m})$ for $3\leq n\leq 9$. 
According to Theorem \ref{th3} we set $m=r+kn^2$ 
where $1<r<n^2$ and $m\neq 0,\pm 1$ is square free.
If $r$ has a common square factor with $n$,
then none of $r+kn^2$ is square-free, we omit $r$. 
If $r$ has no common square factor with $n$ but contains another square factor,
then in our computations it is
replaced by $r+n^2$ or $r+2n^2$ etc. which is already square-free,
but the case will still be represented by $r$ and will cover 
fields $K=\Q(\sqrt[n]{m})$ with square-free $m=r+kn^2$.

For all $m$ we give the integral basis $B$ and discriminant $D$ of $K$.
As far as it is possible we display the index form $I(x_2,\ldots,x_n)$ 
corresponding to
the integral basis and discuss the monogenity of $K$.
Mentioning here the index form equation we always mean the equation
(\ref{iiixxx}).

Stating that the 
index form equation is not solvable modulo $q$ 
in certain cases $m=r+kn^2$ (with a fixed $r$)
we mean that if we let 
$x_2,\ldots,x_n$ and $k$ run through all residue
classes modulo $q$ we never have 
$I(x_2,\ldots,x_n)\equiv \pm 1 \ (\bmod \ q)$. This implies that 
the corresponding fields admit no power integral bases, are not
monogenic.

\vspace{1cm}

\subsection{Pure cubic fields, $K=\Q(\sqrt[3]{m})$}

\vspace{1cm}

In these cases the index form equation is a cubic Thue equation. 

\vspace{1cm}

\noindent
{\bf Case 3.1.} $r=2,3,4,5,6,7$, $m=r+9k$ square-free 
\[
B=\{1,x,x^2\},\;\; D=-27m^2
\]
\[
I(x_2,x_3)= x_2^3-mx_3^3
\]
These fields are obviously monogenic, $(1,0)$ is a solution
of the index form equation.

\vspace{1cm}

\noindent
{\bf Case 3.2.} $r=1$, $m=1+9k$ square-free 
\[
B=\left\{1,x,\frac{1+x+x^2}{3}\right\},\;\; D=-3m^2
\]
\[
I(x_2,x_3)=3x_2^3+3x_2^2x_3+x_2x_3^2-kx_3^3
\]
In this case the index form equation is solvable e.g. for 
$k=27,37$ but not solvable e.g. for $k=10,11,12$.

\vspace{1cm}

\noindent
{\bf Case 3.3.} $r=8$, $m=8+9k$ square-free 
\[
B=\left\{1,x,\frac{1+2x+x^2}{3}\right\},\;\; D=-3m^2
\]
\[
I(x_2,x_3)=3x_2^3+6x_2^2x_3+4x_2x_3^2-kx_3^3
\]
In this case the index form equation is solvable e.g. for 
$k=1,4,12$ but not solvable e.g. for $k=2,3,5,6,7$.

\vspace{1cm}

\subsection{Pure quartic fields, $K=\Q(\sqrt[4]{m})$}

\vspace{1cm}

In these cases the index form is the product of a quadratic and a quartic form.

\vspace{1cm}

\noindent
{\bf Case 4.1.} $r=2,3,6,7,10,11,14,15$ , $m=r+16k$ square-free 
\[
B=\{1,x,x^2,x^3\},\;\; D=-256m^3
\]
\[ 
I(x_2,x_3,x_4)=(x_2^2-mx_4^2)(x_2^4+2mx_2^2x_4^2+m^2x_4^4+4mx_3^4-8mx_2x_4x_3^2)
\]
These fields are monogenic, $(1,0,0)$ is a solution of the index form equation.
This follows also from 
A.Hameed, T.Nakahara, S.M.Husnine and S.Ahmad \cite{nnn}.

\vspace{1cm}

\noindent
{\bf Case 4.2.} $r=1,9$, $m=1+8k$ square-free
\[
B=\left\{1,x,\frac{1+x^2}{2},\frac{1+x+x^2+x^3}{4}\right\},\;\; D=-4m^3
\]
\[
I(x_2,x_3,x_4)=(-x_2x_4-2x_2^2+x_4^2k)\cdot
\]
\[
(x_4^4k^2-2x_4^3x_2k-16x_4^2x_3x_2k+4x_4^2x_2^2k+8x_4^2x_3^2k-16x_4x_3^2x_2k+16x_4x_3^3k+8x_3^4k
\]
\[
+2x_4^2x_2^2-2x_4^2x_3x_2+x_4^2x_3^2-2x_4x_3^2x_2+4x_2^3x_4+2x_4x_3^3+x_3^4+4x_2^4)
\]
If $m=1+16\ell$, that is $k=2\ell$ then the index form equation
is not solvable modulo 2.\\
If $m=9+16\ell$, that is $k=2\ell+1$ then the index form equation has a 
solution for $\ell=4,5$ that is for $m=73,89$ (the solution is (2,1,1)).
For other parameters we conjecture that the minimal index of $K$ is 8.

\vspace{1cm}

\noindent
{\bf Case 4.3.} $r=5,13$, $m=5+8k$ square-free
\[
B=\left\{1,x,\frac{1+x^2}{2},\frac{x+x^3}{2}\right\},\;\; D=-16m^3
\]
\[
I(x_2,x_3,x_4)=
(-x_2x_4-x_2^2+2x_4^2k+x_4^2)\cdot
\]
\[
(16x_4^4k^2+24x_4^4k+16x_4^3x_2k-16x_4^2x_3^2k+16x_4^2x_2^2k-32x_4x_3^2x_2k+8x_3^4k
\]
\[
+9x_4^4+12x_4^3x_2-10x_4^2x_3^2+16x_4^2x_2^2+8x_2^3x_4-20x_4x_3^2x_2+4x_2^4+5x_3^4)
\]
Denote by $f_1$ and $f_2$ the first and second factor of the index form,
respectively. Then we have
\[
f_2-4f_1^2=(8k+5)(2x_2x_4-x_3^2+x_4^2)^2.
\]
If $K$ is monogenic, then for some $x_2,x_3,x_4$ we have $f_1,f_2=\pm 1$,
hence 
\mbox{$f_2-4f_1^2=-3$} or $-5$. This number can only be divisible by $8k+5$ 
for $k=-1$. The field $K=\Q(\sqrt[4]{-3})$ is monogenic, e.g. $(1,1,0)$
is a solution of the index form equation. All other fields of this type
are not monogenic.

\vspace{1cm}

Summarizing the above statements we have
\begin{theorem}
\label{tquartic}
For the following values of $r$ let $m=r+16k$ ($k\in\Z$) be a square-free integer.
The field $K=\Q(\sqrt[4]{m})$ is monogenic for \\
$r=2,3,6,7,10,11,14,15$ and is not monogenic for $r=1,5,13$ with the exception of
$\Q(\sqrt[4]{-3})$ which is monogenic.
\end{theorem}

\begin{remark}
We conjecture that for $m=9+16k$ the only monogenic fields are 
$\Q(\sqrt[4]{73})$, $\Q(\sqrt[4]{89})$.
\end{remark}

\vspace{1cm}

\subsection{Pure quintic fields, $K=\Q(\sqrt[5]{m})$}

\vspace{1cm}

\noindent
{\bf Case 5.1.} $r=2,3,4,5,6,8,9,10,11,12,13,14,15,16,17,19,20,21,22,23$, $m=r+25k$ square-free 
\[
B=\{1,x,x^2,x^3,x^4\},\;\; D=3125m^4.
\]
These fields are obviously monogenic.
The index form is very complicated:
\[
I(x_2,x_3,x_4,x_5)=
-75m^4x_2x_3^2x_4^3x_5^4+45m^4x_2^2x_3x_4^2x_5^5+40m^4x_2x_3^3x_4x_5^5-40m^4x_2x_3x_4^5x_5^3
\]
\[
-75m^2x_2^4x_3^3x_4^2x_5-40m^2x_2^3x_3^5x_4x_5+45m^2x_2^5x_3^2x_4x_5^2+40m^2x_2^5x_3x_4^3x_5+75m^3x_2^3x_3x_4^4x_5^2
\]
\[
+75m^3x_2^2x_3^4x_4x_5^3+50m^3x_2^4x_3x_4x_5^4-200m^3x_2^3x_3^2x_4^2x_5^3+200m^3x_2^2x_3^3x_4^3x_5^2
\]
\[
-45m^3x_2^2x_3^2x_4^5x_5-45m^3x_2x_3^5x_4^2x_5^2-50m^3x_2x_3^4x_4^4x_5-20m^5x_3^2x_4x_5^7+5m^5x_2x_3x_5^8
\]
\[
+35m^5x_3x_4^3x_5^6-15m^5x_2x_4^2x_5^7-5m^4x_2^3x_4x_5^6+20m^4x_2x_4^7x_5^2+25m^4x_2^2x_4^4x_5^4
\]
\[
-25m^4x_3^4x_4^2x_5^4+25m^4x_3^3x_4^4x_5^3-5m^4x_3x_4^8x_5-10m^4x_2^2x_3^2x_5^6+10m^4x_3^2x_4^6x_5^2
\]
\[
-15mx_2^7x_3^2x_5-20mx_2^7x_3x_4^2+5mx_2^8x_4x_5+35mx_2^6x_3^3x_4+20m^2x_2^2x_3^7x_5-5m^2x_2^6x_3x_5^3
\]
\[
+25m^2x_2^4x_3^4x_5^2-25m^2x_2^4x_3^2x_4^4+25m^2x_2^3x_3^4x_4^3-5m^2x_2x_3^8x_4-10m^2x_2^6x_4^2x_5^2
\]
\[
+10m^2x_2^2x_3^6x_4^2-35m^3x_2^3x_4^6x_5+15m^3x_2^2x_3x_4^7-35m^3x_2x_3^6x_5^3+5m^3x_2x_3^3x_4^6+15m^3x_3^7x_4x_5^2
\]
\[
+5m^3x_3^6x_4^3x_5-25m^3x_2^4x_4^3x_5^3-25m^3x_2^3x_3^3x_5^4+x_2^{10}-m^4x_4^{10}-m^2x_3^{10}+x_5^{10}m^6
\]
\[
-11m^5x_4^5x_5^5+11m^4x_3^5x_5^5-11mx_2^5x_3^5+11m^2x_2^5x_4^5-2m^3x_2^5x_5^5+2m^3x_3^5x_4^5
\]

\vspace{1cm}

\noindent
{\bf Case 5.2.} $r=1$, $m=1+25k$ square-free 
\[
B=\left\{1,x,x^2,x^3,\frac{1+x+x^2+x^3+x^4}{5}\right\},\;\; D=125m^4
\]

\vspace{1cm}

\noindent
{\bf Case 5.3.} $r=7$, $m=7+25k$ square-free 
\[             
B=\left\{1,x,x^2,x^3,\frac{1+3x+4x^2+2x^3+x^4}{5}\right\},\;\; D=125m^4
\]
Set $m=7+25k$. For $k=0$ the tuple $(0,-2,-1,2)$ is a solution
of the index form equation,
therefore $K=\Q(\sqrt[5]{7})$ is monogenic. 
(For $k=1$ the $m=32$ is not square free.)

\vspace{1cm}

\noindent
{\bf Case 5.4.} $r=18$, $m=18+25k$ square-free 
\[
B=\left\{1,x,x^2,x^3,\frac{1+2x+4x^2+3x^3+x^4}{5}\right\},\;\; D=125m^4
\]

\vspace{1cm}

\noindent
{\bf Case 5.5.} $r=24$, $m=24+25k$ square-free 
\[
B=\left\{1,x,x^2,x^3,\frac{1+4x+x^2+4x^3+x^4}{5}\right\},\;\; D=125m^4
\]

\vspace{1cm}

\begin{remark}
For the following values of $r$ let $m=r+25k$ ($k\in\Z$) be a square-free integer.
$K=\Q(\sqrt[5]{m})$ is monogenic for 
$r=2,3,4,5,6,8,9,10,11,12,13,14,\\15,16,17,19,20,21,22,23$.
We conjecture that for $r=1,7,18,24$ the fields $\Q(\sqrt[5]{m})$ are not monogenic having minimal index 5 with the exception of $\Q(\sqrt[5]{7})$
which is monogenic.
\end{remark}

\vspace{1cm}

\subsection{Pure sextic fields, $K=\Q(\sqrt[6]{m})$}

\vspace{1cm}

In all these cases the index form is the product of three
factors of degrees  3,6,6, respectively.
We shall denote these factors by $f_1,f_2,f_3$.
These depend on the parameter $m$ and on the variables
$x_2,\ldots,x_6$. 

\vspace{1cm}

\noindent
{\bf Case 6.1.} $r=2,3,6,7,11,14,15,22,23,30,31,34$, $m=r+36k$ square-free 
\[
B=\{1,x,x^2,x^3,x^4,x^5\},\;\; D=6^6 m^5
\]
In this case $K$ is monogenic.
This follows also from the result of 
S.Ahmad, T.Nakahara and S.M.Husnine \cite{n1}
since $m\equiv 2,3 \ (\bmod \ 4)$ and
$m\not \equiv \pm 1 \ (\bmod \ 9)$.
We have 
\[
I(x_2,x_3,x_4,x_5)=(-3mx_2x_4x_6+x_2^3+mx_4^3+m^2x_6^3)
\]
\[
(18m^2x_2^2x_3x_5x_6^2-18m^2x_2x_3^2x_5^2x_6-3m^3x_3^2x_6^4-2m^2x_2^3x_6^3+3m^2x_2^2x_5^4+3m^2x_3^4x_6^2
\]
\[
+2m^2x_3^3x_5^3-3mx_2^4x_5^2-6m^3x_2x_5^2x_6^3+6m^3x_3x_5^3x_6^2-6mx_2^3x_3^2x_6+6mx_2^2x_3^3x_5+m^4x_6^6
\]
\[
-m^3x_5^6-mx_3^6+x_2^6)
\]
\[
(x_2^6+64m^2x_4^6+m^4x_6^6+27m^3x_5^6+27mx_3^6-216m^2x_3^3x_4x_5x_6-72mx_2^3x_3x_4x_5+12m^3x_2x_4x_6^4
\]
\[
+108m^3x_4^2x_5^2x_6^2-108m^3x_4x_5^4x_6+36m^2x_2^2x_4^2x_6^2-96m^2x_2x_4^4x_6+144m^2x_2x_4^3x_5^2
\]
\[
+144m^2x_3^2x_4^3x_6+324m^2x_3^2x_4^2x_5^2-288m^2x_3x_4^4x_5+12mx_2^4x_4x_6+108mx_2^2x_3^2x_4^2
\]
\[
-108mx_2x_3^4x_4-18m^3x_2x_5^2x_6^3+54m^3x_3x_5^3x_6^2-18mx_2^3x_3^2x_6+54mx_2^2x_3^3x_5
\]
\[
-72m^3x_3x_4x_5x_6^3-216m^2x_2x_3x_4x_5^3-54m^2x_2^2x_3x_5x_6^2+162m^2x_2x_3^2x_5^2x_6
\]
\[
-16m^3x_4^3x_6^3-16mx_2^3x_4^3+9m^3x_3^2x_6^4+2m^2x_2^3x_6^3+27m^2x_2^2x_5^4+27m^2x_3^4x_6^2
\]
\[
-54m^2x_3^3x_5^3+9mx_2^4x_5^2)
\]

\vspace{1cm}

\noindent
{\bf Case 6.2.} $r=5,13,21,25,29,33$, $m=r+36k$ square-free
\[
\left\{1,x,x^2,\frac{1+x^3}{2},\frac{x+x^4}{2},\frac{x^2+x^5}{2}\right\},
\;\; D=3^6 m^5
\]
Calculating the index form it is easily seen that
the index form equation is not solvable modulo 2,
hence these fields are not monogenic.
This also follows (in a much more complicated way) from the theorem of
S. Ahmad, T. Nakahara and S. M. Husnine \cite{n2} since all
these $m$ are of the form $m\equiv 1 \ (\bmod \ 4)$ and
$m\not \equiv \pm 1 \ (\bmod \ 9)$.

\vspace{1cm}

\noindent
{\bf Case 6.3.} $r=10,19$, $m=r+36k$ square-free
\[
\left\{1,x,x^2, x^3, \frac{1+x^2+x^4}{3},\frac{x+x^3+x^5}{3}\right\},
\;\; D=2^6 3^2 m^5
\]
Calculating the index form it is easily seen that
the index form equation is not solvable modulo 3,
hence these fields are not monogenic.
This case is not covered by 
S. Ahmad, T. Nakahara and S. M. Husnine \cite{n2}, \cite{n1}  
since $m\equiv 1 \ (\bmod \ 9)$.

\vspace{1cm}

\noindent
{\bf Case 6.4.} $r=26,35$, $m=r+36k$ square-free
\[
\left\{1,x,x^2, x^3,\frac{1+2x^2+x^4}{3},\frac{x+2x^3+x^5}{3}\right\},
\;\; D=2^6 3^2 m^5
\]
This case is not covered by 
S. Ahmad, T. Nakahara and S. M. Husnine \cite{n2}, \cite{n1}  
since $m\equiv -1 \ (\bmod \ 9)$.

If $m=26+36k$ then $4m|(f_2-9f_3)$. If $K$ is monogenic then
for a solution of the index form equation these factors 
are equal to $\pm 1$. The possible values of $f_2-9f_3$
are $\pm 8,\pm 10$, hence the above divisibility can not hold.

If $m=35+36k$ then $4(35+36k)|(f_2-9f_3)$. If $K$ is monogenic then
for a solution of the index form equation these factors 
are equal to $\pm 1$. The possible values of $f_2-9f_3$
are $\pm 8,\pm 10$, hence the above divisibility
can only hold for $k=-1$, that is $m=-1$ which we have excluded.

\vspace{1cm}

\noindent
{\bf Case 6.5.} $r=17$, $m=r+36k$ square-free
\[
\left\{1,x,x^2,\frac{1+x^3}{2}, \frac{4+3x+2x^2+x^4}{6},
\frac{4x+3x^2+2x^3+x^5}{6}\right\},\;\; D= 3^2 m^5
\]
This case is not covered by 
S. Ahmad, T. Nakahara and S. M. Husnine \cite{n2}, \cite{n1}  
since $m\equiv -1 \ (\bmod \ 9)$.
Calculating the index form we can easily see that
the index form equation is not solvable modulo 6.

\vspace{1cm}

\noindent
{\bf Case 6.6.} $m=1$, $m=1+36k$ square-free
\[
\left\{1,x,x^2,\frac{1+x^3}{2},\frac{4+3x+4x^2+x^4}{6},
\frac{3+4x+3x^2+x^3+x^5}{6}\right\},\;\; D= 3^2 m^5
\]
This case is not covered by 
S. Ahmad, T. Nakahara and S. M. Husnine \cite{n2}, \cite{n1}  
since $m\equiv 1 \ (\bmod \ 9)$.
Calculating the index form we can easily see that
the index form equation is not solvable modulo 3.

\vspace{1cm}

Summarizing the above statements we have
\begin{theorem}
\label{tsextic}
For the following values of $r$ let $m=r+36k$ ($k\in\Z$) be a square-free integer.
The field $K=\Q(\sqrt[6]{m})$ is monogenic for \\
$r=2,3,6,7,11,14,15,22,23,30,31,34$ 
and is not monogenic for \\
$r=1,5,10,13,17,19,21,25,26,29,33,35$.
\end{theorem}

\vspace{1cm}

\section{Pure septic fields, $K=\Q(\sqrt[7]{m})$}

\vspace{1cm}

{\bf Case 7.1.} $r=2, 3, 4, 5, 6, 7, 8, 9, 10, 11, 12, 
  13, 14, 15, 16, 17, 20, 21, 22, 23, 24,\\ 25, 26, 27, 28, 29, 32, 
  33, 34, 35, 36, 37, 38, 39, 40, 41, 42, 43, 44, 45, 46, 47$,\\ $m=r+49k$ square-free\\
\[
B=\{1,x,x^2,x^3,x^4,x^5,x^6\},\;\; D=-m^6 7^7
\]
These fields are obviously monogenic.

\vspace{1cm}

{\bf Case 7.2.} $r=18$, $m=18+49k$ square-free\\
\[
B=\{1,x,x^2,x^3,x^4,x^5,\frac{1+2x+4x^2+x^3+2x^4+4x^5+x^6}{7}\},\;\; D=-m^6 7^5
\]

\vspace{1cm}

{\bf Case 7.3.} $r=19$, $m=19+49k$ square-free\\
\[
B=\{1,x,x^2,x^3,x^4,x^5,\frac{1+3x+2x^2+6x^3+4x^4+5x^5+x^6}{7}\},\;\; D=-m^6 7^5
\]

\vspace{1cm}

{\bf Case 7.4.} $r=30$, $m=30+49k$ square-free\\
\[
B=\{1,x,x^2,x^3,x^4,x^5,\frac{1+4x+2x^2+x^3+4x^4+2x^5+x^6}{7}\},\;\; D=-m^6 7^5
\]

\vspace{1cm}

{\bf Case 7.5.} $r=31$, $m=31+49k$ square-free\\
\[
B=\{1,x,x^2,x^3,x^4,x^5,\frac{1+5x+4x^2+6x^3+2x^4+3x^5+x^6}{7}\},\;\; D=-m^6 7^5
\]

\vspace{1cm}

{\bf Case 7.6.} $r=48$, $m=48+49k$ square-free\\
\[
B=\{1,x,x^2,x^3,x^4,x^5,\frac{1+6x+x^2+6x^3+x^4+6x^5+x^6}{7}\},\;\; D=-m^6 7^5
\]

\vspace{1cm}

{\bf Case 7.7.} $r=1$, $m=1+49k$ square-free\\
\[
B=\{1,x,x^2,x^3,x^4,x^5,\frac{1+x+x^2+x^3+x^4+x^5+x^6}{7}\},\;\; D=-m^6 7^5
\]

\vspace{1cm}

\section{Pure octic fields, $K=\Q(\sqrt[8]{m})$}

\vspace{1cm}

In all these cases the index form is the product of three
factors of degrees  4,8,16, respectively.
We shall denote these factors by $f_1,f_2,f_3$.
These depend on the parameter $m$ and on the variables
$x_2,\ldots,x_8$.

\vspace{1cm}

\noindent
{\bf Case 8.1.} $r=2,3,6,7,10,11,14,15,18,19,22,23,26,27,30,31,34,35,38,\\39,42,43,46,
47,50,51,54,55,58,59,62,63$, $m=r+64k$ square-free
\[
\left\{1,x,x^2,x^3,x^4,x^5,x^6,x^7,x^8\right\},\;\; D=-8^8m^7
\]
These fields are obviously monogenic.
This also follows from
A.Hameed, T.Nakahara, S.M.Husnine and S.Ahmad  \cite{nnn}.

\vspace{1cm}

\noindent
{\bf Case 8.2.} $r=1,17,33,49$, $m=r+64k$ square-free
\[
\left\{1,x,x^2,x^3,\frac{1+x^4}{2},\frac{x+x^5}{2},\frac{1+x^2+x^4+x^6}{4},\right.
\]
\[
\left.
\frac{1+x+x^2+x^3+x^4+x^5+x^6+x^7}{8}\right\},\;\; D=-2^{10}m^7
\]       
These fields are not monogenic by the theorem of 
A.Hameed and T.Nakahara \cite{n8}.                           
We conjecture that in these fields the minimal index is 128.

\vspace{1cm}

\noindent
{\bf Case 8.3.} $r=5,13,21,29,37,45,53,61$, these cases can be included by $m=5+8k$,
square-free
\[
\left\{1,x,x^2,x^3,\frac{1+x^4}{2},\frac{x+x^5}{2},\frac{x^2+x^6}{2},
\frac{x^3+x^7}{2}\right\},\;\; D=-2^{16}m^7
\]

Calculating and factorizing $f_3-16f_2^2$ we find that it is divisible
by $m$. If there existed a power integral basis then for 
a solution of the index form equation
we would have $f_1,f_2,f_3=\pm 1$, hence
$f_3-16f_2^2=\pm 1-16$ is either $-15$ or $-17$.
The possible divisors are $\pm 3,\pm 5,\pm 15,\pm 17$ but only $m=-3,5$
is of type $m=5+8k$. 

For $m=-3$ then element $(-1,-1,0,1,1,0,-1)$ has index one, hence 
$K=\Q(\sqrt[8]{-3})$ is monogenic.

For $m=5$ the least index we found in $K=\Q(\sqrt[8]{5})$ was 16.

Note that A.Hameed, T.Nakahara, S.M.Husnine and S.Ahmad  
\cite{n8} assert that these fields are not monogenic, they 
certainly did not involve $K=\Q(\sqrt[8]{-3})$.

\vspace{1cm}

\noindent
{\bf Case 8.4.} $r=9,25,41,57$, these cases can be included by $m=9+16k$,
square-free
\[
\left\{1,x,x^2,x^3,\frac{1+x^4}{2},\frac{x+x^5}{2},\frac{1+x^2+x^4+x^6}{4},
\frac{x+x^3+x^5+x^7}{4}\right\},\;\; D=-2^{12}m^7
\]
Calculating and factorizing $f_2-4f_1^2$ we find that it is divisible
by $m$. If there existed a power integral basis then for 
a solution of the index form equation we would have $f_1,f_2,f_3=\pm 1$, hence
$f_2-4f_1^2=\pm 1-4$ is either $-3$ or $-5$.
The possible divisors are $\pm 3,\pm 5$ but none of them 
is of type $m=9+16k$. 
Therefore these fields are not monogenic.
This also follows from the theorem of 
A.Hameed, T.Nakahara, S.M.Husnine and S.Ahmad \cite{n8}.

\vspace{1cm}

Summarizing the above statements we have
\begin{theorem}
\label{toctic}
For the following values of $r$ let $m=r+64k$ ($k\in\Z$) be a square-free integer.
The field $K=\Q(\sqrt[8]{m})$ is monogenic for \\
$r=2,3,6,7,10,11,14,15,18,19,22,23,26,27,30,31,34,35,38,39,42,43,46$,\\
$47,50,51,54,55,58,59,62,63$ 
and is not monogenic for \\
$r=1,5,9,13,17,21,25,29,33,37,41,45,49,53,57,61$, $m\neq 5$,
with the exception of $K=\Q(\sqrt[8]{-3})$ which is monogenic.
\end{theorem}

\begin{remark}
We conjecture that the minimal index of $K=\Q(\sqrt[8]{5})$ is 16.
Octic fields of this type will be considered in a forecoming paper.
\end{remark}

\vspace{1cm}

\section{Pure nonic fields, $K=\Q(\sqrt[9]{m})$}

\vspace{1cm}

\noindent
{\bf Case 9.1.} $r=2,3,4,5,6,7,11,12,13,14,15,16,20,21,22,23,24,25,29,30,\\31,32,33,34,
38,39,40,41,42,43,47,48,49,50,51,52,56,57,58,59,60,61,65,\\66,67,68,69,70,74,
75,76,77,78,79$, $m=r+81k$, square-free
\[
\left\{1,x,x^2,x^3,x^4,x^5,x^6,x^7,x^8\right\},\;\; D=3^{18}m^8
\]
These fields are obviously monogenic.

\vspace{1cm}

\noindent
{\bf Case 9.2.} $r=1,28,55$, $m=r+81k$, square-free
\[
\left\{1,x,x^2,x^3,x^4,x^5,\frac{1+x^3+x^6}{3},\frac{x+x^4+x^7}{3},
\right.
\]
\[
\left.
\frac{1+x+x^2+x^3+x^4+x^5+x^6+x^7+x^8}{9}\right\},\;\; D=3^{10}m^8
\]

\vspace{1cm}

\noindent
{\bf Case 9.3.} $r=8,17,35,44,62,71 $, $m=r+81k$, square-free
\[
\left\{1,x,x^2,x^3,x^4,x^5,\frac{1+2x^3+x^6}{3},\frac{x+2x^4+x^7}{3},
\frac{x^2+2x^5+x^8}{3}\right\},\;\; D=3^{12}m^8
\]

\vspace{1cm}

\noindent
{\bf Case 9.4.} $r=10,19,37,46,64,73$, $m=r+81k$, square-free
\[
\left\{1,x,x^2,x^3,x^4,x^5,\frac{1+x^3+x^6}{3},\frac{x+x^4+x^7}{3},
\frac{x^2+x^5+x^8}{3}\right\},\;\; D=3^{12}m^8
\]

\vspace{1cm}

\noindent
{\bf Case 9.5.} $r=26,53,80$, $m=r+81k$, square-free
\[
\left\{1,x,x^2,x^3,x^4,x^5,\frac{1+2x^3+x^6}{3},\frac{x+2x^4+x^7}{3},\right.
\]
\[
\left.
\frac{1+2x+x^2+8x^3+7x^4+8x^5+x^6+2x^7+x^8}{9}\right\}, D=3^{10}m^8
\]

\vspace{1cm}

\section{Computational remarks}
In all our calculations we used Maple \cite{maple}
and most of our programs executed a couple of seconds or a few minutes on 
an average laptop. For $n=4,6,8$ we needed a very careful calculation
of the factors of the index forms, which may take extremely long otherwise.

The tests corresponding to Theorem \ref{th3} took also a few minutes
for $n=3,4,5,6,8$. For $n=9$ it executed 5 hours.
For $n=7$ we executed our Maple program on a supercomputer 
with nodes having 24 CPU-s.
The running time on one node was 10 hours per remainder.
We had 48 remainders and the program was running on 10 nodes
parallelly.


\begin{thebibliography}{10}

\normalsize
\baselineskip=17pt

\bibitem{n3}
A.Hameed and T.Nakahara,
{\em Integral bases and relative monogenity of pure octic fields},
Bull. Math. Soc. Sci. Math. R\'epub. Soc. Roum., Nouv. S\'er. 
{\bf 58(106)}(2015) No.4, 419--433.

\bibitem{n8}
A.Hameed and T.Nakahara,
{\em On pure octic fields related to a problem of Hasse},
manuscript.

\bibitem{n2}
S. Ahmad, T. Nakahara and S. M. Husnine, 
{\em Non-monogenesis of a family of pure sextic fields}, 
Arch. Sci. (Geneva) {\bf 65}(2012), No. 7, 42–-49.

\bibitem{n1}
S.Ahmad, T.Nakahara and S.M.Husnine,
{\em Power integral bases for certain pure sextic fields},
Int. J. Number Theory {\bf 10}(2014), No. 8, 2257--2265.

\bibitem{nnn}
A.Hameed, T.Nakahara, S.M.Husnine and S.Ahmad, 
{\em On existence of canonical number system in certain classes of pure algebraic number fields}, 
J. Prime Research in Mathematics, {\bf 7}(2011), 19--24.

\bibitem{maple}
B.W.Char, K.O.Geddes, G.H.Gonnet, M.B.Monagan, S.M.Watt (eds.)
{\em MAPLE, Reference Manual}, Watcom Publications, Waterloo, Canada, 1988.


\bibitem{fun}
T.Funakura,
{\em On integral bases of pure quartic fields},
Math. J. Okayama Univ. {\bf 26}(1984), 27--41.

\bibitem{gaal}
Diophantine equations and power integral bases,
Boston, Birkh\"auser, 2002.

\bibitem{gr}
I.Ga\'al and L.Remete,
{\em Binomial Thue equations and power integral bases in pure quartic fields}
JP Journal of Algebra Number Theory Appl. {\bf 32}(2014), No. 1, 49--61.

\bibitem{grsz}
I.Ga\'al, L.Remete and T.Szab\'o,
{\em Calculating power integral bases by using relative power integral bases}
Functiones et Approximatio, to appear.

\bibitem{gsz}
I.Ga\'al and T.Szab\'o,
{\em A note on the minimal indices of pure cubic fields},
JP Journal of Algebra Number Theory Appl. {\bf 19}(2010), No. 2, 129--139.


\bibitem{fad}
L. El Fadil, 
{\em Computation of a power integral basis of a pure cubic number field},
Int. J. Contemp. Math. Sci. {\bf 2}(2007), No. 13--16, 601--606.

\bibitem{sqy}
B.K.Spearman, Y.Qiduan and J.Yoo,
{\em Minimal indices of pure cubic fields},
Arch. Math. {\bf 106}(2016), No.1, 35--40.

\bibitem{sp}
B.K.Spearman and K.S.Williams, 
{\em An explicit integral basis for a pure cubic field},
Far East J. Math. Sci. {\bf 6}(1998), No. 1, 1--14.


\end{thebibliography}
\end{document}